\documentclass{amsart}
\def\N{{\mathbb{N}}}
\def\Z{{\mathbb{Z}}}
\def\R{{\mathbb{R}}}

\def\Image{\operatorname{Im}}

\def\Cl{\operatorname{Cl}}

\newtheorem{Theorem}{Theorem}[section]

\newtheorem{Lemma}[Theorem]{Lemma}
\newtheorem{Proposition}[Theorem]{Proposition}

\theoremstyle{definition}
\newtheorem{Definition}[Theorem]{Definition}

\theoremstyle{remark}

\begin{document}
\sloppy
\title{On equivariant homeomorphisms of boundaries of CAT(0) groups}
\author{Tetsuya Hosaka} 
\address{Department of Mathematics, Faculty of Education, 
Utsunomiya University, Utsunomiya, 321-8505, Japan}
\date{October 3, 2010}
\email{hosaka@cc.utsunomiya-u.ac.jp}
\thanks{
Partly supported by the Grant-in-Aid for Young Scientists (B), 
The Ministry of Education, Culture, Sports, Science and Technology, Japan.
(No.\ 21740037).}
\begin{abstract}
In this paper, 
we investigate an equivariant homeomorphism of 
the boundaries $\partial X$ and $\partial Y$ of 
two proper CAT(0) spaces $X$ and $Y$ on which a CAT(0) group $G$ acts geometrically.
We provide a sufficient condition to obtain a $G$-equivariant homeomorphism of 
the two boundaries $\partial X$ and $\partial Y$ 
as a continuous extension of the quasi-isometry 
$\phi:Gx_0\rightarrow Gy_0$ defined by $\phi(gx_0)=gy_0$, 
where $x_0\in X$ and $y_0\in Y$.
\end{abstract}
\maketitle
%
%%%%%%%%%%%%%
% Section 1 %
%%%%%%%%%%%%%
\section{Introduction}

In this paper, 
we investigate an equivariant homeomorphism of 
the boundaries of two proper CAT(0) spaces on which a CAT(0) group acts geometrically 
as a continuous extension of a quasi-isometry of the two CAT(0) spaces.

Definitions and details of CAT(0) spaces and their boundaries 
are found in \cite{BH} and \cite{GH}.
A {\it geometric} action on a CAT(0) space 
is an action by isometries which is proper (\cite[p.131]{BH}) and cocompact.
We note that every CAT(0) space on which some group acts 
geometrically is a proper space (\cite[p.132]{BH}).
A group $G$ is called a {\it CAT(0) group}, 
if $G$ acts geometrically on some CAT(0) space $X$.

It is well-known that 
if a Gromov hyperbolic group $G$ acts 
geometrically on a negatively curved space $X$, then 
the natural map $G\rightarrow X$ $(g\mapsto gx_0)$ extends continuously to 
an equivariant homeomorphism of the boundaries of $G$ and $X$.
Also if a Gromov hyperbolic group $G$ acts 
geometrically on negatively curved spaces $X$ and $Y$, then 
the boundaries of $X$ and $Y$ are $G$-equivariant homeomorphic.
Indeed 
the natural map $Gx_0 \rightarrow Gy_0$ $(gx_0\mapsto gy_0)$ extends continuously to 
a $G$-equivariant homeomorphism of the boundaries of $X$ and $Y$.
The boundaries of Gromov hyperbolic groups are quasi-isometric invariant 
(cf.\ \cite{BH}, \cite{CP}, \cite{GH}, \cite{Gr}, \cite{Gr0}).

Here in \cite{Gr0}, Gromov asked whether 
the boundaries of two CAT(0) spaces $X$ and $Y$ are $G$-equivariant homeomorphic 
whenever a CAT(0) group $G$ acts geometrically on the two CAT(0) spaces $X$ and $Y$.
In \cite{BR}, 
P.~L.~Bowers and K.~Ruane have constructed an example that 
the natural quasi-isometry $Gx_0 \rightarrow Gy_0$ $(gx_0\mapsto gy_0)$ 
does not extend continuously to 
any map between the boundaries $\partial X$ and $\partial Y$ of $X$ and $Y$.
Also S.~Yamagata \cite{Y} has constructed a similar example 
using a right-angled Coxeter group and its Davis complex.
Moreover, there is a research by C.~Croke and B.~Kleiner \cite{CK0} on 
an equivariant homeomorphism of the boundaries $\partial X$ and $\partial Y$.

Also, C.~Croke and B.~Kleiner \cite{CK} have constructed a CAT(0) group $G$ 
which acts geometrically on two CAT(0) spaces $X$ and $Y$ 
whose boundaries are not homeomorphic, 
and J.~Wilson \cite{W} has proved 
that this CAT(0) group has uncountably many boundaries.
Recently, C.~Mooney \cite{Moo1} has showed that 
the knot group $G$ of any connected sum of two non-trivial torus knots 
has uncountably many CAT(0) boundaries.

Also, it has been observed by M.~Bestvina \cite{Bes} 
that all the boundaries of a given CAT(0) group are shape equivalent, 
and he has asked the question whether 
all the boundaries of a given CAT(0) group are cell-like equivalent.
This question is an open problem and there are some resent research (cf.\ \cite{AGW}, \cite{Moo2}).

The purpose of this paper is 
to provide a sufficient condition to obtain a $G$-equivariant homeomorphism between 
the two boundaries $\partial X$ and $\partial Y$ of 
two CAT(0) spaces $X$ and $Y$ on which a CAT(0) group $G$ acts geometrically 
as a continuous extension of the natural quasi-isometry 
$Gx_0\rightarrow Gy_0$ $(gx_0\mapsto gy_0)$, 
where $x_0\in X$ and $y_0\in Y$.

Now we recall the example of Bowers and Ruane in \cite{BR}.
Let $G= F_2\times \Z$ and $X=Y=T\times \R$, 
where $F_2$ is the rank 2 free group generated by $\{a,b\}$ and 
$T$ is the Cayley graph of $F_2$ with respect to the generating set $\{a,b\}$.
Then we define the action ``$\cdot$'' of the group $G$ on the CAT(0) space $X$ by 
\begin{align*}
&(a,0)\cdot(t,r)=(a\cdot t,r), \\
&(b,0)\cdot(t,r)=(b\cdot t,r), \\
&(1,1)\cdot(t,r)=(t,r+1),
\end{align*}
for each $(t,r)\in T\times \R=X$, 
and also define the action ``$*$'' of the group $G$ on the CAT(0) space $Y$ by 
\begin{align*}
&(a,0)*(t,r)=(a\cdot t,r), \\
&(b,0)*(t,r)=(b\cdot t,r+2), \\ 
&(1,1)*(t,r)=(t,r+1),
\end{align*}
for each $(t,r)\in T\times \R=Y$.
Then the group $G$ acts geometrically on the two CAT(0) spaces $X$ and $Y$, 
and the quasi-isometry $g\cdot x_0 \mapsto g*y_0$ 
(where $x_0=(1,0)\in X$ and $y_0=(1,0)\in Y$) 
does not extend continuously to 
any map from $\partial X$ to $\partial Y$.
Indeed for $g_i=a^ib^i\in F_2$, 
$\{g_i^\infty\,|\,i\in\N\}\rightarrow a^\infty$ as $i\rightarrow\infty$ in $\partial T$,
\begin{align*}
&\lim_{n\rightarrow \infty}(g_i^n,0)\cdot x_0=[g_i^\infty,0], \\
&\lim_{n\rightarrow \infty}(a^n,0)\cdot x_0=[a^\infty,0],
\end{align*}
in $X\cup\partial X$, and 
\begin{align*}
&\lim_{n\rightarrow \infty}(g_i^n,0)*y_0=[g_i^\infty,\frac{\pi}{4}], \\
&\lim_{n\rightarrow \infty}(a^n,0)*y_0=[a^\infty,0],
\end{align*}
in $Y\cup\partial Y$.
Hence any map from $\partial X$ to $\partial Y$ obtained 
as a continuously extension of the quasi-isometry 
$G\cdot x_0\rightarrow G*y_0$ $(g\cdot x_0\rightarrow g*y_0)$ 
must send $[g_i^\infty,0]$ to $[g_i^\infty,\frac{\pi}{4}]$ and 
fix $[a^\infty,0]$.
However, this is incompatible with continuously at $[a^\infty,0]$, 
because $[g_i^\infty,0]\rightarrow [a^\infty,0]$ as $i\rightarrow \infty$ 
(\cite[p.187]{BR}).

Here in this example, we note that 
\begin{enumerate}
\item[(a)] the point $a^i\cdot x_0$ is in the geodesic segment from $x_0$ to $g^i\cdot x_0$ in $X$, 
i.e., $a^i\cdot x_0\in [x_0,g^i\cdot x_0]$ in $X$ for any $i\in\N$ and 
\item[(b)] the distance between the point $a^i*y_0$ and the geodesic segment from $y_0$ to $g^i*y_0$ 
is unbounded for $i\in\N$ in $Y$, i.e., 
there does {\it not} exist a constant $M>0$ such that 
$d(a^i*y_0, [y_0,g^i*y_0])\le M$ for any $i\in\N$ in $Y$.
\end{enumerate}

Based on this observation, we consider a condition.

We suppose that a group $G$ acts geometrically on two CAT(0) spaces $X$ and $Y$.
Let $x_0\in X$ and $y_0\in Y$.
Then we define the condition $(*)$ as follows:
\begin{enumerate}
\item[$(*)$] 
There exist constants $N>0$ and $M>0$ such that $GB(x_0,N)=X$, $GB(y_0,M)=Y$ and 
for any $g,a\in G$, 
if $[x_0,gx_0] \cap B(ax_0,N)\neq\emptyset$ in $X$ then 
$[y_0,gy_0] \cap B(ay_0,M)\neq\emptyset$ in $Y$.
\end{enumerate}

In this paper, we prove the following theorem.

\begin{Theorem}
Suppose that a group $G$ acts geometrically on two CAT(0) spaces $X$ and $Y$.
Let $x_0\in X$ and $y_0\in Y$.
If the condition~$(*)$ holds, then 
there exists a $G$-equivariant homeomorphism of 
the boundaries $\partial X$ and $\partial Y$ 
as a continuous extension of the quasi-isometry 
$\phi:Gx_0\rightarrow Gy_0$ defined by $\phi(gx_0)=gy_0$.
\end{Theorem}

%%%%%%%%%%%%%
% Section 2 %
%%%%%%%%%%%%%
\section{CAT(0) spaces and their boundaries}

Details of CAT(0) spaces and their boundaries 
are found in \cite{ABN}, \cite{BH}, \cite{GO}, \cite{GH} and \cite{Sw}.

A proper geodesic space $(X,d_X)$ is called a {\it CAT(0) space}, 
if the ``CAT(0)-inequality'' holds 
for all geodesic triangles $\triangle$ and 
for all choices of two points $x$ and 
$y$ in $\triangle$.
Here the ``CAT(0)-inequality'' is defined as follows: 
Let $\triangle$ be a geodesic triangle in $X$.
A {\it comparison triangle} for $\triangle$ is 
a geodesic triangle $\triangle'$ in the Euclidean plain $\R^2$
with same edge lengths as $\triangle$.
Choose two points $x$ and $y$ in $\triangle$. 
Let $x'$ and $y'$ denote the corresponding points in $\triangle'$.
Then the inequality $$d_X(x,y) \le d_{\R^2}(x',y')$$ 
is called the {\it CAT(0)-inequality}, 
where $d_{\R^2}$ is the natural metric on $\R^2$.

Every proper CAT(0) space can be compactified by 
adding its ``boundary''.
Let $(X,d_X)$ be a proper CAT(0) space, 
and let $\mathcal{R}$ be the set of all geodesic rays in $X$.
We define an equivalence relation $\sim$ in $\mathcal{R}$ as follows:
For geodesic rays $\xi,\zeta:[0,\infty)\rightarrow X$, 
$$ \xi\sim\zeta \iff 
\Image \xi\subset B(\Image\zeta,N) \ \text{for some $N\ge0$},$$
where $B(A,N):=\{x\in X\,|\, d_X(x,A)\le N\}$ for $A\subset X$.
Then the {\it boundary} $\partial X$ of $X$ is 
defined as 
$$ \partial X = \mathcal{R}/\sim. $$
For each geodesic ray $\xi\in\mathcal{R}$, 
the equivalence class of $\xi$ is denoted by $\xi(\infty)$.

It is known that for each $\alpha\in\partial X$ and each $x_0 \in X$, 
there exists a unique geodesic ray $\xi_{\alpha}:[0,\infty)\rightarrow X$ 
such that $\xi_{\alpha}(0)=x_0$ and $\xi_{\alpha}(\infty)=\alpha$.
Thus we can identify the boundary $\partial X$ of $X$ as 
the set of all geodesic rays $\xi$ with $\xi(0)=x_0$.

Let $(X,d_X)$ be a proper CAT(0) space and let $x_0\in X$.
We define a topology on $X\cup \partial X$ as follows: 
\begin{enumerate}
\item[(1)] $X$ is an open subspace of $X\cup\partial X$. 
\item[(2)] Let $\alpha\in\partial X$ and let $\xi_{\alpha}$ be the geodesic ray such that 
$\xi_{\alpha}(0)=x_0$ and $\xi_{\alpha}(\infty)=\alpha$.
For $r>0$ and $\epsilon>0$, we define 
$$ U_{X \cup \partial X}(\alpha;r,\epsilon)=
\{ x \in X \cup \partial X \,|\,
x \not\in B(x_0,r),\ d_X(\xi_{\alpha}(r),\xi_{x}(r))<\epsilon \}, $$
where $\xi_{x}:[0,d_X(x_0,x)]\rightarrow X$ is 
the geodesic (segment or ray) from $x_0$ to $x$.
Let $\epsilon_0>0$ be a constant.
Then the set 
$$\{U_{X\cup\partial X}(\alpha;r,\epsilon_0)\,|\,r>0 \} $$
is a neighborhood basis for $\alpha$ in $X\cup\partial X$.
\end{enumerate}
Here it is known that the topology on $X\cup \partial X$ 
is not dependent on the basepoint $x_0 \in X$ and 
$X\cup\partial X$ is a metrizable compactification of $X$.

Also 
for $\alpha\in\partial X$ and the geodesic ray $\xi_{\alpha}$ with 
$\xi_{\alpha}(0)=x_0$ and $\xi_{\alpha}(\infty)=\alpha$ 
and for $r>0$ and $\epsilon>0$, 
we define 
$$ U'_{X \cup \partial X}(\alpha;r,\epsilon)=
\{ x \in X \cup \partial X \,|\,
x \not\in B(x_0,r),\ d_X(\xi_{\alpha}(r),\Image\xi_{x})<\epsilon \}, $$
where $\xi_{x}:[0,d(x_0,x)]\rightarrow X$ is 
the geodesic (segment or ray) from $x_0$ to $x$.
Let $\epsilon_0>0$ be a constant.
Then the set 
$$\{U'_{X\cup\partial X}(\alpha;r,\epsilon_0)\,|\,r>0 \} $$
is also a neighborhood basis for $\alpha$ in $X\cup\partial X$ 
(cf.\ \cite[Lemma~4.2]{Ho000}).

Suppose that a group $G$ acts on a proper CAT(0) space $X$ by isometries.
For each element $g \in G$ and 
each geodesic ray $\xi:[0,\infty)\rightarrow X$, 
a map $g \xi:[0,\infty)\rightarrow X$ 
defined by $(g\xi)(t):=g(\xi(t))$ is also a geodesic ray.
For two geodesic rays $\xi$ and $\xi'$, 
if $\xi(\infty)=\xi'(\infty)$ 
then $g\xi(\infty)=g\xi'(\infty)$. 
Thus $g$ induces a homeomorphism of $\partial X$, 
and $G$ acts on $\partial X$ by homeomorphisms.
Here we note that 
if a sequence $\{x_i\,|\,i\in \N\}\subset X$ 
converges to $\alpha\in\partial X$ in $X\cup\partial X$, then 
for any $g\in G$, 
the sequence $\{gx_i\,|\,i\in \N\}\subset X$ 
converges to $g\alpha\in\partial X$ in $X\cup\partial X$.

\begin{Definition}%%%%%%%%%%%%%%%%%%%%%%%%%%%%%%% "Cauchy sequence"
Let $(X,d_X)$ be a proper CAT(0) space and 
let $\{x_i\,|\,i\in\N\}\subset X$ be an unbounded sequence in $X$.
In this paper, we say that 
the sequence $\{x_i\,|\,i\in\N\}$ is a {\it Cauchy sequence in $X\cup\partial X$}, 
if there exists $\epsilon_0>0$ such that for any $r>0$, 
there is a number $i_0\in\N$ as 
$$ x_i \in U_{X\cup\partial X}(x_{i_0};r,\epsilon_0) $$
for any $i\ge i_0$.
Here 
$$ U_{X\cup\partial X}(x_{i_0};r,\epsilon_0)= 
\{ x \in X \,|\, x \not\in B(x_0,r),\ d_X(\xi_{x_{i_0}}(r),\xi_{x}(r))<\epsilon \}, $$
where $\xi_{z}$ is the geodesic segment from $x_0$ to $z$ in $X$.
\end{Definition}

We show the following lemma which is used later.

\begin{Lemma}\label{Lem1}
Let $(X,d_X)$ be a proper CAT(0) space and 
let $\{x_i\,|\,i\in\N\}\subset X$ be an unbounded sequence in $X$.
Then the sequence $\{x_i\,|\,i\in\N\}$ is a Cauchy sequence in $X\cup\partial X$ defined above 
if and only if 
the sequence $\{x_i\,|\,i\in\N\}$ converges to some point $\alpha\in\partial X$ in $X\cup\partial X$.
\end{Lemma}

\begin{proof}
We first show that 
if the sequence $\{x_i\,|\,i\in\N\}$ converges to some point $\alpha\in\partial X$ in $X\cup\partial X$, 
then $\{x_i\,|\,i\in\N\}$ is a Cauchy sequence in $X\cup\partial X$ defined above.

Suppose that $\{x_i\,|\,i\in\N\}$ converges to $\alpha\in\partial X$ in $X\cup\partial X$.
Let $\epsilon_0>0$.
Since the set 
$$\{U_{X\cup\partial X}(\alpha;r,\frac{\epsilon_0}{2})\,|\,r>0\}$$ 
is a neighborhood basis for $\alpha$ in $X\cup\partial X$, for each $r>0$, 
there exists a number $i_0\in\N$ such that 
$$x_i\in U_{X\cup\partial X}(\alpha;r,\frac{\epsilon_0}{2})$$ 
for any $i\ge i_0$.
Then for any $i\ge i_0$, 
\begin{align*}
d_X(\xi_{x_{i_0}}(r),\xi_{x}(r)) &\le 
d_X(\xi_{x_{i_0}}(r),\xi_{\alpha}(r)) + d_X(\xi_{\alpha}(r),\xi_{x}(r)) \\
&\le \frac{\epsilon_0}{2}+\frac{\epsilon_0}{2} \\
&=\epsilon_0.
\end{align*}
Hence 
$x_i \in U_{X\cup\partial X}(x_{i_0};r,\epsilon_0)$ for any $i\ge i_0$.
Thus the sequence $\{x_i\,|\,i\in\N\}$ is a Cauchy sequence in $X\cup\partial X$.

Next, we show that 
if $\{x_i\,|\,i\in\N\}$ is a Cauchy sequence in $X\cup\partial X$ defined above, 
then 
$\{x_i\,|\,i\in\N\}$ converges to some point $\alpha\in\partial X$ in $X\cup\partial X$.

Suppose that $\{x_i\,|\,i\in\N\}$ is a Cauchy sequence in $X\cup\partial X$.
Since the set $\{x_i\,|\,i\in\N\}$ is unbounded in $X$, 
there exists a limit point $\alpha\in \Cl\{x_i\,|\,i\in\N\} \cap \partial X$.
Here there exists a subsequence $\{x_{i_j}\,|\,j\in \N\}\subset \{x_i\,|\,i\in\N\}$ 
which converges to $\alpha$ in $X\cup\partial X$.

Then we show that 
the sequence $\{x_i\,|\,i\in\N\}$ converges to the point $\alpha\in\partial X$ in $X\cup\partial X$.

Since $\{x_i\,|\,i\in\N\}$ is a Cauchy sequence in $X\cup\partial X$, 
there exists $\epsilon_0>0$ such that for any $r>0$, 
there is a number $i_0\in\N$ as 
$x_i \in U_{X\cup\partial X}(x_{i_0};r,\epsilon_0)$ for any $i\ge i_0$, i.e., 
$d_X(\xi_{x_{i_0}}(r),\xi_{x_i}(r))\le \epsilon_0$ for any $i\ge i_0$.
Also since the subsequence $\{x_{i_j}\,|\,j\in \N\}$ converges to $\alpha$ in $X\cup\partial X$, 
there exists $i_{j_0}\ge i_0$ such that 
$x_{i_{j_0}}\in U_{X\cup\partial X}(\alpha;r,1)$, i.e., 
$d_X(\xi_{x_{i_{j_0}}}(r),\xi_{\alpha}(r))\le 1$.
Then for any $i\ge i_{j_0}$, 
\begin{align*}
d_X(\xi_{x_i}(r),\xi_{\alpha}(r)) &\le 
d_X(\xi_{x_i}(r),\xi_{x_{i_0}}(r)) + d_X(\xi_{x_{i_0}}(r),\xi_{x_{i_{j_0}}}(r)) 
+ d_X(\xi_{x_{i_{j_0}}}(r),\xi_{\alpha}(r)) \\
&\le \epsilon_0 + \epsilon_0 + 1 \\
&= 2\epsilon_0+1,
\end{align*}
since $i\ge i_{j_0}\ge i_0$.
Hence for any $r>0$, 
there exists a number $i_{j_0}\in\N$ such that for any $i\ge i_{j_0}$, 
$$x_i\in U_{X\cup\partial X}(\alpha;r,2\epsilon_0+1),$$ 
where $2\epsilon_0+1$ is a constant.
Thus the sequence $\{x_i\,|\,i\in\N\}$ converges 
to the point $\alpha\in\partial X$ in $X\cup\partial X$.
\end{proof}

%%%%%%%%%%%%%
% Section 3 %
%%%%%%%%%%%%%
\section{Proof of the main theorem}

We suppose that a group $G$ acts geometrically on two CAT(0) spaces $(X,d_X)$ and $(Y,d_Y)$.
Let $x_0\in X$ and $y_0\in Y$.
Now we suppose that the condition $(*)$ holds; that is, 
\begin{enumerate}
\item[$(*)$] 
there exist constants $N>0$ and $M>0$ such that $GB(x_0,N)=X$, $GB(y_0,M)=Y$ and 
for any $g,a\in G$, 
if $[x_0,gx_0] \cap B(ax_0,N)\neq\emptyset$ in $X$ then 
$[y_0,gy_0] \cap B(ay_0,M)\neq\emptyset$ in $Y$.
\end{enumerate}

Our goal is to show that 
the quasi-isometry $\phi:Gx_0\rightarrow Gy_0$ defined by $\phi(gx_0)=gy_0$ 
continuously extends to a $G$-equivariant homeomorphism 
of the boundaries $\partial X$ and $\partial Y$.

Since the map $\phi:Gx_0\rightarrow Gy_0$ defined by $\phi(gx_0)=gy_0$ 
is a quasi-isometry (cf.\ \cite[p.138]{BH}, \cite{Gr}, \cite{Gr0}), 
there exist constants $\lambda>0$ and $C>0$ such that 
$$ \frac{1}{\lambda} d_Y(gy_0,hy_0)-C \le d_X(gx_0,hx_0) \le \lambda d_Y(gy_0,hy_0)+C $$
for any $g,h\in G$.

We first show the following.

\begin{Proposition}\label{Prop2-1}
Let $\{g_i\}\subset G$ be a sequence.
If $\{g_ix_0\}\subset X$ is a Cauchy sequence in $X\cup \partial X$ defined in Section~2, 
then 
$\{g_iy_0\}\subset Y$ is also a Cauchy sequence in $Y\cup \partial Y$.
\end{Proposition}

\begin{proof}
Let $\{g_i\}\subset G$.
Suppose that $\{g_ix_0\}\subset X$ is a Cauchy sequence in $X\cup \partial X$.

To prove that $\{g_iy_0\}\subset Y$ is a Cauchy sequence in $Y\cup \partial Y$, 
we show that 
there exists $M'>0$ such that for any $R>0$, 
there is $i_0\in \N$ as 
$$ g_iy_0 \in U_{Y\cup\partial Y}(g_{i_0}y_0;R,M') $$
for any $i\ge i_0$.

Let $M'=\lambda(2N+1)+2M+C$ and let $R>0$.

Since $\{g_ix_0\}\subset X$ is a Cauchy sequence in $X\cup \partial X$, 
for $r=\lambda(R+C+M)+N$, 
there exists $i_0\in \N$ such that 
$$ g_ix_0 \in U_{X\cup\partial X}(g_{i_0}x_0;r,1) $$
for any $i\ge i_0$.

Then 
$$d_X(x_0,g_{i_0}x_0)\ge r, \ d_X(x_0,g_ix_0)\ge r, \ \text{and}\ 
d_X(\xi_{g_{i_0}x_0}(r),\xi_{g_ix_0}(r)) \le 1,$$
where $\xi_{g_{i_0}x_0}$ is the geodesic from $x_0$ to $g_{i_0}x_0$ and 
$\xi_{g_ix_0}$ is the geodesic from $x_0$ to $g_ix_0$ in $X$.

Since $GB(x_0,N)=X$, 
there exist $a,b\in G$ such that 
$d_X(ax_0,\xi_{g_{i_0}x_0}(r))\le N$ and $d_X(bx_0,\xi_{g_ix_0}(r))\le N$.
Then 
$$ [x_0,g_{i_0}x_0]\cap B(ax_0,N)\neq \emptyset \ \text{and} \ 
[x_0,g_ix_0]\cap B(bx_0,N)\neq \emptyset.$$
Hence by the condition $(*)$, 
$$ [y_0,g_{i_0}y_0]\cap B(ay_0,M)\neq \emptyset \ \text{and} \ 
[y_0,g_iy_0]\cap B(by_0,M)\neq \emptyset.$$
Thus 
$$\xi_{g_{i_0}y_0}(r'_0)\in [y_0,g_{i_0}y_0]\cap B(ay_0,M) \ \text{and} \ 
\xi_{g_iy_0}(r')\in [y_0,g_iy_0]\cap B(by_0,M)$$ 
for some $r'_0>0$ and $r'>0$.

To obtain that for any $i\ge i_0$, 
$$ g_iy_0 \in U_{Y\cup\partial Y}(g_{i_0}y_0;R,M'), $$
we show that 
$$ r'_0 \ge R, \ r'\ge R \ \text{and}\ 
d_Y(\xi_{g_{i_0}y_0}(r'_0),\xi_{g_iy_0}(r'))\le M'.$$

First, 
\begin{align*}
r'_0 &= d_Y(y_0,\xi_{g_{i_0}y_0}(r'_0)) \\
&\ge d_Y(y_0,ay_0)-M \\
&\ge \frac{1}{\lambda}d_X(x_0,ax_0)-C-M \\
&\ge \frac{1}{\lambda}(r-N)-C-M \\
&=R,
\end{align*}
because $d_X(x_0,ax_0)\ge r-N$ and $r=\lambda(R+C+M)+N$.

By the same argument,
\begin{align*}
r' &= d_Y(y_0,\xi_{g_iy_0}(r')) \\
&\ge d_Y(y_0,by_0)-M \\
&\ge \frac{1}{\lambda}d_X(x_0,bx_0)-C-M \\
&\ge \frac{1}{\lambda}(r-N)-C-M \\
&=R,
\end{align*}
because $d_X(x_0,bx_0)\ge r-N$ and $r=\lambda(R+C+M)+N$.

Also, 
\begin{align*}
d_Y(\xi_{g_{i_0}y_0}(r'_0),\xi_{g_iy_0}(r')) 
&\le d_Y(ay_0,by_0)+2M \\
&\le (\lambda d_X(ax_0,bx_0)+C)+2M \\
&\le \lambda(d_X(\xi_{g_{i_0}x_0}(r),\xi_{g_ix_0}(r))+2N)+C+2M \\
&\le \lambda(1+2N)+C+2M \\
&= M',
\end{align*}
because $d_X(\xi_{g_{i_0}x_0}(r),\xi_{g_ix_0}(r))\le 1$ and 
$M'=\lambda(2N+1)+2M+C$.

Thus 
$$ r'_0 \ge R, \ r'\ge R \ \text{and}\ 
d_Y(\xi_{g_{i_0}y_0}(r'_0),\xi_{g_iy_0}(r'))\le M'.$$
Hence 
$$ 
d_Y(\xi_{g_{i_0}y_0}(R),\xi_{g_iy_0}(R))\le 
d_Y(\xi_{g_{i_0}y_0}(r'_0),\xi_{g_iy_0}(r')) \le M',$$
since $Y$ is a CAT(0) space.
Also we obtain that 
$$d_Y(y_0,g_{i_0}y_0)\ge R \ \text{and} \ d_Y(y_0,g_iy_0)\ge R,$$
because $r'_0 \ge R$ and $r'\ge R$.

Thus 
$$ g_iy_0 \in U_{Y\cup\partial Y}(g_{i_0}y_0;R,M') $$
for any $i\ge i_0$.
Hence we obtain that 
$\{g_iy_0\}\subset Y$ is a Cauchy sequence in $Y\cup \partial Y$.
\end{proof}

Then we define a map 
$\bar{\phi}:\partial X\rightarrow \partial Y$ 
as a continuous extension of the quasi-isometry 
$\phi:Gx_0\rightarrow Gy_0$ defined by $\phi(gx_0)=gy_0$ as follows:
For each $\alpha\in \partial X$, 
there exists a sequence $\{g_ix_0\}\subset Gx_0\subset X$ 
which converges to $\alpha$ in $X\cup\partial X$.
Then the sequence $\{g_ix_0\}\subset X$ is a Cauchy sequence in $X\cup \partial X$ 
by Lemma~\ref{Lem1}.
By Proposition~\ref{Prop2-1}, 
the sequence $\{g_iy_0\}\subset Y$ is also a Cauchy sequence in $Y\cup\partial Y$.
Hence by Lemma~\ref{Lem1}, 
the sequence $\{g_iy_0\}\subset Y$ converges to 
some point $\bar{\alpha}\in\partial Y$ in $Y\cup\partial Y$.
Then we define $\bar{\phi}(\alpha)=\bar{\alpha}$.

\begin{Proposition}\label{Prop2-2}
The map $\bar{\phi}:\partial X\rightarrow \partial Y$ is well-defined.
\end{Proposition}

\begin{proof}
Let $\alpha\in \partial X$ and 
let $\{g_ix_0\},\{h_ix_0\} \subset Gx_0\subset X$ be two sequences 
which converge to $\alpha$ in $X\cup\partial X$.
As the argument above, by Lemma~\ref{Lem1} and Proposition~\ref{Prop2-1}, 
the sequence $\{g_iy_0\}\subset Y$ converges to 
some point $\bar{\alpha}\in\partial Y$ and 
the sequence $\{h_iy_0\}\subset Y$ converges to 
some point $\bar{\beta}\in\partial Y$ in $Y\cup\partial Y$.
Then we show that $\bar{\alpha}=\bar{\beta}$.

Here we can consider a sequence $\{\tilde{g}_jx_0\,|\,j\in\N\}\subset Gx_0\subset X$ 
such that 
$$\{\tilde{g}_jx_0\,|\,j\in\N\}=\{g_ix_0\,|\,i\in\N\}\cup \{h_ix_0\,|\,i\in\N\}$$
and the sequence $\{\tilde{g}_jx_0\}$ converges to $\alpha$ in $X\cup\partial X$.
Then the sequence $\{\tilde{g}_jx_0\}$ is a Cauchy sequence in $X\cup\partial X$ and 
the sequence $\{\tilde{g}_jy_0\}$ is also in $Y\cup\partial Y$ by Proposition~\ref{Prop2-1}.
Hence the sequence $\{\tilde{g}_jy_0\}$ converges to 
some point $\bar{\gamma}\in\partial Y$ in $Y\cup\partial Y$.
Here we note that 
the two sequences $\{g_iy_0\}$ and $\{h_iy_0\}$ are subsequences of $\{\tilde{g}_jy_0\}$.
Hence we obtain that $\bar{\alpha}=\bar{\beta}=\bar{\gamma}$.

Thus the map $\bar{\phi}:\partial X\rightarrow \partial Y$ defined as above 
is well-defined.
\end{proof}

Next, we show the following.

\begin{Proposition}\label{Prop2-3}
The map $\bar{\phi}:\partial X\rightarrow \partial Y$ is surjective.
\end{Proposition}

\begin{proof}
Let $\bar{\alpha}\in \partial Y$.
There exists a sequence $\{g_iy_0\}\subset Gy_0 \subset Y$ 
which converges to $\bar{\alpha}$ in $Y\cup \partial Y$.
Then we consider the set $\{g_ix_0\,|\,i\in \N\}$ 
which is an unbounded subset of $X$.
Here 
$$ \Cl\{g_ix_0\,|\,i\in \N\} \cap \partial X \neq \emptyset,$$
and there exists a subsequence $\{g_{i_j}x_0\,|\,j\in\N\}\subset \{g_ix_0\}$ 
which converges to some point $\alpha \in \partial X$.
Then the sequence $\{g_{i_j}y_0\}$ converges to $\bar{\alpha}$ in $Y\cup \partial Y$, 
because $\{g_{i_j}y_0\}$ is a subsequence of the sequence $\{g_iy_0\}$ 
which converges to $\bar{\alpha}$ in $Y\cup \partial Y$.
Hence $\bar{\phi}(\alpha)=\bar{\alpha}$ 
by the definition of the map $\bar{\phi}$.
Thus the map $\bar{\phi}:\partial X\rightarrow \partial Y$ is surjective.
\end{proof}

Here we provide a lemma.

\begin{Lemma}\label{Lem2-4}
For any $\tilde{N}\ge N$, 
there exists $\tilde{M}>0$ such that $GB(y_0,\tilde{M})=Y$ and 
for any $g,a\in G$, 
if $[x_0,gx_0] \cap B(ax_0,\tilde{N})\neq\emptyset$ in $X$ then 
$[y_0,gy_0] \cap B(ay_0,\tilde{M})\neq\emptyset$ in $Y$.
\end{Lemma}

\begin{proof}
For $\tilde{N}\ge N$, 
we put $\tilde{M}=\lambda(N+\tilde{N})+C+M$.

Let $g,a\in G$ 
as $[x_0,gx_0] \cap B(ax_0,\tilde{N})\neq\emptyset$ in $X$.
Then there exists a point $x_1\in [x_0,gx_0] \cap B(ax_0,\tilde{N})$.
Since $GB(x_0,N)=X$, 
there exists $a'\in G$ such that $x_1 \in B(a'x_0,N)$.
Then $x_1 \in [x_0,gx_0]\cap B(a'x_0,N)$ and 
$[x_0,gx_0]\cap B(a'x_0,N)\neq\emptyset$ in $X$.
By the condition~$(*)$, 
$[y_0,gy_0]\cap B(a'y_0,M)\neq\emptyset$ in $Y$.
Hence $d_Y(a'y_0,[y_0,gy_0])\le M$.
Here we note that 
\begin{align*}
d_Y(a'y_0,ay_0) &\le \lambda d_X(a'x_0,ax_0)+C \\
&\le \lambda(d_X(a'x_0,x_1)+d_X(x_1,ax_0))+C \\
&\le \lambda(N+\tilde{N})+C.
\end{align*}
Hence 
\begin{align*}
d_Y(ay_0,[y_0,gy_0]) &\le d_Y(ay_0,a'y_0)+d_Y(a'y_0,[y_0,gy_0]) \\
&\le \lambda(N+\tilde{N})+C+M \\
&=\tilde{M}.
\end{align*}
Thus we obtain that 
$[y_0,gy_0] \cap B(ay_0,\tilde{M})\neq\emptyset$ in $Y$.
\end{proof}

Let $\tilde{N}=2N$.
By Lemma~\ref{Lem2-4}, 
there exists $\tilde{M}>0$ such that $GB(y_0,\tilde{M})=Y$ and 
for any $g,a\in G$, 
if $[x_0,gx_0] \cap B(ax_0,\tilde{N})\neq\emptyset$ in $X$ then 
$[y_0,gy_0] \cap B(ay_0,\tilde{M})\neq\emptyset$ in $Y$.

Here we show the following technical lemma.

\begin{Lemma}\label{Lem2-5}
Let $\alpha\in\partial X$ and let $\xi_{\alpha}:[0,\infty)\rightarrow X$ 
be the geodesic ray in $X$ such that 
$\xi_{\alpha}(0)=x_0$ and $\xi_{\alpha}(\infty)=\alpha$.
Let $\{g_ix_0\}\subset Gx_0 \subset X$ be a sequence 
which converges to $\alpha$ in $X\cup\partial X$ 
such that $d_X(g_ix_0,\xi_{\alpha}(i))\le N$ for any $i\in\N$ 
(since $GB(x_0,N)=X$, we can take such a sequence).
Then 
\begin{enumerate}
\item[$(1)$] $d_X(g_ix_0,[x_0,g_jx_0]) \le \tilde{N}$ for any $i,j\in \N$ with $i<j$,
\item[$(2)$] $d_Y(g_iy_0,[y_0,g_jy_0]) \le \tilde{M}$ for any $i,j\in \N$ with $i<j$,
\item[$(3)$] $d_Y(g_iy_0,\Image \xi_{\bar{\alpha}}) \le \tilde{M}+1$ for any $i\in\N$,
\item[$(4)$] $d_X(g_ix_0,g_{i+1}x_0)\le 2N+1$ for any $i\in\N$,
\item[$(5)$] $d_Y(g_iy_0,g_{i+1}y_0)\le \lambda(2N+1)+C$ for any $i\in\N$, and 
\item[$(6)$] $\Image\xi_{\bar{\alpha}} \subset 
\bigcup\{B(g_iy_0,3(\tilde{M}+1)+\lambda(2N+1)+C)\,|\,i\in\N \}$.
\end{enumerate}
Here $\bar{\alpha}=\bar{\phi}(\alpha)$ and 
$\xi_{\bar{\alpha}}:[0,\infty)\rightarrow Y$ is the geodesic ray in $Y$ such that 
$\xi_{\bar{\alpha}}(0)=y_0$ and $\xi_{\bar{\alpha}}(\infty)=\bar{\alpha}$.
\end{Lemma}

\begin{proof}
(1) 
For any $i,j\in \N$ with $i<j$, 
\begin{align*}
d_X(g_ix_0,[x_0,g_jx_0]) 
&\le d_X(g_ix_0,\xi_{\alpha}(i))+d_X(\xi_{\alpha}(i),[x_0,g_jx_0]) \\
&\le N+N =2N \\
&=\tilde{N},
\end{align*}
where we obtain the inequality $d_X(\xi_{\alpha}(i),[x_0,g_jx_0])\le N$, 
since $d_X(g_jx_0,\xi_{\alpha}(j))\le N$, $i<j$ and $X$ is a CAT(0) space.

(2) 
By Lemma~\ref{Lem2-4} and the definition of $\tilde{M}$, 
we obtain that 
$d_Y(g_iy_0,[y_0,g_jy_0]) \le \tilde{M}$ for any $i,j\in \N$ with $i<j$ from (1).

(3) 
We note that 
the sequence $\{g_iy_0\}$ converges to $\bar{\alpha}$ 
by the definition of the map $\bar{\phi}:\partial X\rightarrow \partial Y$.

Let $i\in\N$ and let $R=d_Y(y_0,g_iy_0)$.
Since the sequence $\{g_jy_0\}$ converges to $\bar{\alpha}$, 
there exists $j_0\in N$ such that 
$$d_Y(\xi_{\bar{\alpha}}(R),\xi_{g_jy_0}(R))<1$$
for any $j\ge j_0$, 
because the set 
$$\{U_{Y\cup\partial Y}(\bar{\alpha};r,1)\,|\,r>0\}$$ 
defined in Section~2 
is a neighborhood basis for $\bar{\alpha}$ in $Y\cup\partial Y$.

Let $j\in \N$ with $j>i$ and $j>j_0$.
Since $i<j$, we obtain that $d_Y(g_iy_0,[y_0,g_jy_0]) \le \tilde{M}$ by (2).
Hence there exists $r>0$ such that 
$d_Y(g_iy_0,\xi_{g_jy_0}(r)) \le \tilde{M}$.
Here we note that $r \le R$ by \cite[Lemma~4.1]{Ho000}
and we can obtain that 
$$d_Y(\xi_{\bar{\alpha}}(r),\xi_{g_jy_0}(r))<d_Y(\xi_{\bar{\alpha}}(R),\xi_{g_jy_0}(R))<1,$$
since $Y$ is a CAT(0) space.
Then 
\begin{align*}
d_Y(g_iy_0,\Image \xi_{\bar{\alpha}}) 
&\le d_Y(g_iy_0,\xi_{g_jy_0}(r))+d_Y(\xi_{g_jy_0}(r),\Image \xi_{\bar{\alpha}}) \\
&< \tilde{M}+1.
\end{align*}
Hence
$d_Y(g_iy_0,\Image \xi_{\bar{\alpha}}) \le \tilde{M}+1$ for any $i\in\N$.

(4) 
We obtain that $d_X(g_ix_0,g_{i+1}x_0)\le 2N+1$ for any $i\in\N$, 
because 
\begin{align*}
d_X(g_ix_0,g_{i+1}x_0) &\le 
d_X(g_ix_0,\xi_{\alpha}(i))+d_X(\xi_{\alpha}(i),\xi_{\alpha}(i+1))+d_X(\xi_{\alpha}(i+1),g_{i+1}x_0) \\
&\le N+1+N \\
&=2N+1,
\end{align*}
since $d_X(g_ix_0,\xi_{\alpha}(i))\le N$ for any $i\in\N$ 
by the definition of the sequence $\{g_ix_0\}$.

(5) 
Since the map $\phi:Gx_0\rightarrow Gy_0 \ (gx_0\mapsto gy_0)$ is a quasi-isometry, 
we obtain that 
$d_Y(g_iy_0,g_{i+1}y_0)\le \lambda(2N+1)+C$ for any $i\in\N$ by (4).

(6) 
For each $i\in\N$, there exists $r_i>0$ such that 
$d_Y(g_iy_0,\xi_{\bar{\alpha}}(r_i))\le \tilde{M}+1$ by (3).
Then by (5), 
\begin{align*}
d_Y(\xi_{\bar{\alpha}}(r_i),\xi_{\bar{\alpha}}(r_{i+1})) &\le 
d_Y(\xi_{\bar{\alpha}}(r_i),g_iy_0)+d_Y(g_iy_0,g_{i+1}y_0)+
d_Y(g_{i+1}y_0,\xi_{\bar{\alpha}}(r_{i+1})) \\
&\le (\tilde{M}+1)+(\lambda(2N+1)+C)+(\tilde{M}+1) \\
&= 2(\tilde{M}+1)+\lambda(2N+1)+C.
\end{align*}
Hence we obtain that 
$$ \Image\xi_{\bar{\alpha}} \subset 
\bigcup\{B(g_iy_0,3(\tilde{M}+1)+\lambda(2N+1)+C)\,|\,i\in\N \}.$$
\end{proof}

Now we show the following.

\begin{Proposition}\label{Prop2-6}
The map $\bar{\phi}:\partial X\rightarrow \partial Y$ is injective.
\end{Proposition}

\begin{proof}
Let $\alpha,\alpha'\in\partial X$, and 
let $\xi_{\alpha}:[0,\infty)\rightarrow X$ and $\xi_{\alpha'}:[0,\infty)\rightarrow X$ 
be the geodesic rays in $X$ such that 
$\xi_{\alpha}(0)=\xi_{\alpha'}(0)=x_0$, 
$\xi_{\alpha}(\infty)=\alpha$ and $\xi_{\alpha'}(\infty)=\alpha'$.
Let 
$\{g_ix_0\},\{g'_ix_0\}\subset Gx_0 \subset X$ be sequences such that 
$d_X(g_ix_0,\xi_{\alpha}(i))\le N$ and $d_X(g'_ix_0,\xi_{\alpha'}(i))\le N$.
Then 
the sequence $\{g_ix_0\}$ converges to $\alpha$ and 
the sequence $\{g'_ix_0\}$ converges to $\alpha'$ in $X\cup\partial X$.

Let $\bar{\alpha}=\bar{\phi}(\alpha)$ and $\bar{{\alpha'}}=\bar{\phi}(\alpha')$.
Also let $\xi_{\bar{\alpha}}:[0,\infty)\rightarrow Y$ and 
$\xi_{\bar{{\alpha'}}}:[0,\infty)\rightarrow Y$ 
be the geodesic rays in $Y$ such that 
$\xi_{\bar{\alpha}}(0)=\xi_{\bar{{\alpha'}}}(0)=y_0$, 
$\xi_{\bar{\alpha}}(\infty)=\bar{\alpha}$ and 
$\xi_{\bar{{\alpha'}}}(\infty)=\bar{{\alpha'}}$.

Then by Lemma~\ref{Lem2-5}, 
\begin{enumerate}
\item[$(1)$] $d_X(g_ix_0,[x_0,g_jx_0]) \le \tilde{N}$ for any $i,j\in \N$ with $i<j$,
\item[$(2)$] $d_Y(g_iy_0,[y_0,g_jy_0]) \le \tilde{M}$ for any $i,j\in \N$ with $i<j$,
\item[$(3)$] $d_Y(g_iy_0,\Image \xi_{\bar{\alpha}}) \le \tilde{M}+1$ for any $i\in\N$,
\item[$(4)$] $d_X(g_ix_0,g_{i+1}x_0)\le 2N+1$ for any $i\in\N$,
\item[$(5)$] $d_Y(g_iy_0,g_{i+1}y_0)\le \lambda(2N+1)+C$ for any $i\in\N$,
\item[$(6)$] $\Image\xi_{\bar{\alpha}} \subset 
\bigcup\{B(g_iy_0,3(\tilde{M}+1)+\lambda(2N+1)+C)\,|\,i\in\N \}$, 
\end{enumerate}
and 
\begin{enumerate}
\item[$(1')$] $d_X(g'_ix_0,[x_0,g'_jx_0]) \le \tilde{N}$ for any $i,j\in \N$ with $i<j$,
\item[$(2')$] $d_Y(g'_iy_0,[y_0,g'_jy_0]) \le \tilde{M}$ for any $i,j\in \N$ with $i<j$,
\item[$(3')$] $d_Y(g'_iy_0,\Image \xi_{\bar{{\alpha'}}}) \le \tilde{M}+1$ for any $i\in\N$,
\item[$(4')$] $d_X(g'_ix_0,g'_{i+1}x_0)\le 2N+1$ for any $i\in\N$,
\item[$(5')$] $d_Y(g'_iy_0,g'_{i+1}y_0)\le \lambda(2N+1)+C$ for any $i\in\N$,
\item[$(6')$] $\Image\xi_{\bar{{\alpha'}}} \subset 
\bigcup\{B(g'_iy_0,3(\tilde{M}+1)+\lambda(2N+1)+C)\,|\,i\in\N \}$.
\end{enumerate}

To prove that 
the map $\bar{\phi}:\partial X\rightarrow \partial Y$ is injective, 
we show that if $\alpha\neq\alpha'$ then $\bar{\alpha}\neq\bar{{\alpha'}}$.

We suppose that $\alpha\neq\alpha'$.
Then the geodesic rays $\xi_{\alpha}$ and $\xi_{\alpha'}$ are not asymptotic.
Hence for any $t>0$, 
there exists $r_0>0$ such that 
$d_X(\xi_{\alpha}(r_0),\Image \xi_{\alpha'})>t$.
Then for $i_0\in\N$ with $i_0\ge r_0$,
\begin{align*}
d_X(g_{i_0}x_0,\Image \xi_{\alpha'}) 
&\ge 
d_X(\xi_{\alpha}(i_0),\Image \xi_{\alpha'})-d_X(g_{i_0}x_0,\xi_{\alpha}(i_0)) \\
&\ge 
d_X(\xi_{\alpha}(r_0),\Image \xi_{\alpha'})-N \\
&> t-N
\end{align*}
Since $d_X(g'_jx_0,\Image \xi_{\alpha'}) \le N$ for any $j\in \N$, 
we obtain that 
$d_X(g_{i_0}x_0,g'_jx_0) > t-2N$ for any $j\in \N$.
Hence for any $j\in \N$, 
\begin{align*}
d_Y(g_{i_0}y_0,g'_jy_0) &\ge \frac{1}{\lambda}d_X(g_{i_0}x_0,g'_jx_0)-C \\
&> \frac{1}{\lambda}(t-2N)-C.
\end{align*}
Here by $(6')$, 
$$ \Image\xi_{\bar{{\alpha'}}} \subset 
\bigcup\{B(g'_jy_0,3(\tilde{M}+1)+\lambda(2N+1)+C)\,|\,j\in\N \}.$$
Let $j_0\in\N$ such that 
$$d_Y(g_{i_0}y_0,g'_{j_0}y_0)=\min\{d_Y(g_{i_0}y_0,g'_jy_0)\,|\,j\in\N \}.$$
Then 
\begin{align*}
d_Y(g_{i_0}y_0,\Image \xi_{\bar{{\alpha'}}}) &\ge 
\min\{d_Y(g_{i_0}y_0,g'_jy_0)\,|\,j\in\N \}-(3(\tilde{M}+1)+\lambda(2N+1)+C) \\
&= d_Y(g_{i_0}y_0,g'_{j_0}y_0)-(3(\tilde{M}+1)+\lambda(2N+1)+C) \\
&> (\frac{1}{\lambda}(t-2N)-C)-(3(\tilde{M}+1)+\lambda(2N+1)+C),
\end{align*}
since 
$d_Y(g_{i_0}y_0,g'_jy_0)>\frac{1}{\lambda}(t-2N)-C$ for any $j\in \N$ 
by the argument above.

Thus for any $t>0$, 
there exists $i_0\in \N$ such that 
$$ d_Y(g_{i_0}y_0,\Image \xi_{\bar{{\alpha'}}}) > 
(\frac{1}{\lambda}(t-2N)-C)-(3(\tilde{M}+1)+\lambda(2N+1)+C).$$
Here by (3), there exists $R_0>0$ such that 
$$ d_Y(g_{i_0}y_0,\xi_{\bar{\alpha}}(R_0)) \le \tilde{M}+1.$$
Then 
\begin{align*}
d_Y(\xi_{\bar{\alpha}}(R_0),\Image \xi_{\bar{{\alpha'}}}) 
&\ge d_Y(g_{i_0}y_0,\Image \xi_{\bar{{\alpha'}}})-d_Y(g_{i_0}y_0,\xi_{\bar{\alpha}}(R_0)) \\
&> (\frac{1}{\lambda}(t-2N)-C)-(3(\tilde{M}+1)+\lambda(2N+1)+C)-(\tilde{M}+1) \\
&= (\frac{1}{\lambda}(t-2N)-C)-(4(\tilde{M}+1)+\lambda(2N+1)+C).
\end{align*}
Since $t>0$ is an arbitrary large number, 
the two geodesic rays $\xi_{\bar{\alpha}}$ and $\xi_{\bar{{\alpha'}}}$ are not asymptotic 
and $\bar{\alpha}\neq\bar{{\alpha'}}$.

Therefore, 
the map $\bar{\phi}:\partial X\rightarrow \partial Y$ is injective.
\end{proof}

From Propositions~\ref{Prop2-3} and \ref{Prop2-6}, 
we obtain that 
the map $\bar{\phi}:\partial X\rightarrow \partial Y$ is bijective.

We show the following.

\begin{Proposition}
The map $\bar{\phi}:\partial X\rightarrow \partial Y$ is continuous.
\end{Proposition}

\begin{proof}
Let $\alpha\in\partial X$ and let $\bar{\alpha}=\bar{\phi}(\alpha)$.
We put $\bar{c}=\lambda(2N+3)+C+2(\tilde{M}+1)$ which is a constant.

To prove that 
the map $\bar{\phi}:\partial X\rightarrow \partial Y$ is continuous at the point $\alpha\in\partial X$, 
we show that for any $\bar{r}>0$, there exists $r>0$ such that 
if $\beta\in U_{X\cup\partial X}(\alpha;r,1)$ then 
$\bar{\beta}\in U'_{Y\cup\partial Y}(\bar{\alpha};\bar{r},\bar{c})$ 
where $\bar{\beta}=\bar{\phi}(\beta)$, 
because 
$\{U_{X\cup\partial X}(\alpha;r,1)\,|\,r>0\}$ and 
$\{U'_{Y\cup\partial Y}(\bar{\alpha};\bar{r},\bar{c})\,|\,\bar{r}>0\}$ are 
neighborhood basis for $\alpha$ and $\bar{\alpha}$ in 
$\partial X$ and $\partial Y$, respectively.

For $\bar{r}>0$, 
we take $r=\lambda(\bar{r}+C+\tilde{M}+1)+N+1$.

Let $\beta\in U_{X\cup\partial X}(\alpha;r,1)$ 
and let $\bar{\beta}=\bar{\phi}(\beta)$.

By Lemma~\ref{Lem2-5}, 
there exists a sequence $\{g_i\}\subset G$ such that 
\begin{enumerate}
\item[$(1)$] the sequence $\{g_ix_0\}\subset X$ converges to $\alpha$ in $X\cup\partial X$,
\item[$(2)$] $d_X(g_ix_0,\Image \xi_{\alpha}) \le N$ for any $i\in\N$,
\item[$(3)$] $\Image\xi_{\alpha} \subset \bigcup\{B(g_ix_0,N+1)\,|\,i\in\N \}$,
\item[$(4)$] the sequence $\{g_iy_0\}\subset Y$ converges to $\bar{\alpha}$ in $Y\cup\partial Y$, and 
\item[$(5)$] $d_Y(g_iy_0,\Image \xi_{\bar{\alpha}}) \le \tilde{M}+1$ for any $i\in\N$.
\end{enumerate}
Here since $d_X(g_ix_0,\xi_{\alpha}(i))\le N$ and 
$d_X(\xi_{\alpha}(i),\xi_{\alpha}(i+1))=1$ for any $i\in \N$ in Lemma~\ref{Lem2-5}, 
we can obtain the statement $(3)$ above.
Also, 
there exists a sequence $\{h_j\}\subset G$ such that 
\begin{enumerate}
\item[$(1')$] the sequence $\{h_jx_0\}\subset X$ converges to $\beta$ in $X\cup\partial X$,
\item[$(2')$] $d_X(h_jx_0,\Image \xi_{\beta}) \le N$ for any $j\in\N$,
\item[$(3')$] $\Image\xi_{\beta} \subset \bigcup\{B(h_jx_0,N+1)\,|\,j\in\N \}$,
\item[$(4')$] the sequence $\{h_jy_0\}\subset Y$ converges to $\bar{\beta}$ in $Y\cup\partial Y$, and 
\item[$(5')$] $d_Y(h_jy_0,\Image \xi_{\bar{\beta}}) \le \tilde{M}+1$ for any $j\in\N$.
\end{enumerate}

Since $\beta\in U_{X\cup\partial X}(\alpha;r,1)$, 
$$d_X(\xi_{\alpha}(r),\xi_{\beta}(r))<1.$$
By $(3)$ and $(3')$, 
there exist $i_0\in\N$ and $j_0\in\N$ such that 
$$ d_X(g_{i_0}x_0,\xi_{\alpha}(r))\le N+1 \ \text{and} \ 
d_X(h_{j_0}x_0,\xi_{\beta}(r))\le N+1.$$
Also by $(5)$ and $(5')$, 
there exist $\tilde{r}>0$ and $\tilde{{r'}}>0$ such that 
$$ d_Y(g_{i_0}y_0,\xi_{\bar{\alpha}}(\tilde{r}))\le \tilde{M}+1 \ \text{and} \ 
d_Y(h_{j_0}y_0,\xi_{\bar{\beta}}(\tilde{{r'}}))\le \tilde{M}+1.$$

Then 
\begin{align*}
d_X(g_{i_0}x_0,h_{j_0}x_0) &\le 
d_X(g_{i_0}x_0,\xi_{\alpha}(r))+d_X(\xi_{\alpha}(r),\xi_{\beta}(r))
+d_X(\xi_{\beta}(r),h_{j_0}x_0) \\
&\le (N+1)+1+(N+1) \\
&=2N+3.
\end{align*}
Hence 
\begin{align*}
d_Y(g_{i_0}y_0,h_{j_0}y_0) &\le \lambda d_X(g_{i_0}x_0,h_{j_0}x_0)+C \\
&\le \lambda(2N+3)+C.
\end{align*}
Then 
\begin{align*}
d_Y(\xi_{\bar{\alpha}}(\tilde{r}),\xi_{\bar{\beta}}(\tilde{{r'}})) &\le 
d_Y(\xi_{\bar{\alpha}}(\tilde{r}),g_{i_0}y_0)+
d_Y(g_{i_0}y_0,h_{j_0}y_0)+
d_Y(h_{j_0}y_0,\xi_{\bar{\beta}}(\tilde{{r'}})) \\
&\le (\tilde{M}+1)+(\lambda(2N+3)+C)+(\tilde{M}+1) \\
&=\lambda(2N+3)+C+2(\tilde{M}+1) \\
&=\bar{c}.
\end{align*}
Also, 
\begin{align*}
\tilde{r} &= d_Y(y_0,\xi_{\bar{\alpha}}(\tilde{r})) \\
&\ge d_Y(y_0,g_{i_0}y_0)-d_Y(g_{i_0}y_0,\xi_{\bar{\alpha}}(\tilde{r})) \\
&\ge d_Y(y_0,g_{i_0}y_0)-(\tilde{M}+1) \\
&\ge \frac{1}{\lambda}d_X(x_0,g_{i_0}x_0)-C-(\tilde{M}+1) \\
&\ge \frac{1}{\lambda}(d_X(x_0,\xi_{\alpha}(r))-d_X(g_{i_0}x_0,\xi_{\alpha}(r)))-C-(\tilde{M}+1) \\
&\ge \frac{1}{\lambda}(r-(N+1))-C-(\tilde{M}+1) \\
&=\bar{r},
\end{align*}
since $r=\lambda(\bar{r}+C+\tilde{M}+1)+N+1$.
Thus we obtain that 
\begin{align*}
d_Y(\xi_{\bar{\alpha}}(\bar{r}),\Image\xi_{\bar{\beta}}) &\le 
d_Y(\xi_{\bar{\alpha}}(\tilde{r}),\Image\xi_{\bar{\beta}}) \\
&\le d_Y(\xi_{\bar{\alpha}}(\tilde{r}),\xi_{\bar{\beta}}(\tilde{{r'}})) \\
&\le \bar{c}.
\end{align*}
Hence 
$\bar{\beta}\in U'_{Y\cup\partial Y}(\bar{\alpha};\bar{r},\bar{c})$.

Thus the map $\bar{\phi}:\partial X\rightarrow \partial Y$ is continuous.
\end{proof}

Finally, we show the following.

\begin{Theorem}
The map $\bar{\phi}:\partial X\rightarrow \partial Y$ is a $G$-equivariant homeomorphism.
\end{Theorem}

\begin{proof}
By the argument above, 
the map $\bar{\phi}:\partial X\rightarrow \partial Y$ is well-defined, bijective and continuous.

From the definition and the well-definedness of $\bar{\phi}$, 
we obtain that 
the map $\bar{\phi}:\partial X\rightarrow \partial Y$ is $G$-equivariant.
Indeed for any $\alpha\in\partial X$ and $g\in G$, 
if $\{g_ix_0\}\subset Gx_0\subset X$ is a sequence which converges to $\alpha$ in $X\cup\partial X$, 
then $\bar{\phi}(\alpha)$ is the point of $\partial Y$ 
to which the sequence $\{g_iy_0\}\subset Gy_0\subset Y$ converges in $Y\cup\partial Y$.
Then $\{gg_ix_0\}\subset Gx_0\subset X$ is the sequence which converges to $g\alpha$ in $X\cup\partial X$ 
and $\bar{\phi}(g\alpha)$ is the point of $\partial Y$ 
to which the sequence $\{gg_iy_0\}\subset Gy_0\subset Y$ converges in $Y\cup\partial Y$.
Here we note that 
the sequence $\{gg_iy_0\}\subset Gy_0\subset Y$ converges to $g\bar{\phi}(\alpha)$ in $Y\cup\partial Y$ 
by the definition of the action of $G$ on $\partial Y$.
Hence $\bar{\phi}(g\alpha)=g\bar{\phi}(\alpha)$ 
for any $\alpha\in\partial X$ and $g\in G$ and 
the map $\bar{\phi}:\partial X\rightarrow \partial Y$ is $G$-equivariant.

Also, 
the map $\bar{\phi}:\partial X\rightarrow \partial Y$ is closed, 
since $\partial X$ and $\partial Y$ are compact and metrizable.

Therefore, we obtain that 
the map $\bar{\phi}:\partial X\rightarrow \partial Y$ is 
a $G$-equivariant homeomorphism.
\end{proof}

%%%%%%%%%%%%%
% Section 4 %
%%%%%%%%%%%%%
\section{Remark}

The author thinks that 
there is a possibility that 
the main theorem, the condition $(*)$ 
and some arguments in this paper 
can be used to investigate boundaries of CAT(0) groups and interesting open problems on 
\begin{enumerate}
\item[(1)] (equivariant) rigidity of boundaries of CAT(0) groups; 
\item[(2)] (equivariant) rigidity of boundaries of Coxeter groups; 
\item[(3)] (equivariant) rigidity of boundaries of Davis complexes of Coxeter groups; 
\item[(4)] (equivariant) rigidity of boundaries of CAT(0) spaces 
on which Coxeter groups act geometrically by reflections;
\item[(5)] (equivariant) rigidity of boundaries of CAT(0) spaces 
on which right-angled Coxeter groups act geometrically by reflections;
\item[(6)] (equivariant) rigidity of boundaries of CAT(0) cubical complexes 
on which CAT(0) groups act geometrically, 
\end{enumerate}
etc.

Here we can find some recent research on CAT(0) groups and their boundaries 
in \cite{CK}, \cite{GO}, \cite{Ha}, \cite{Ho}, \cite{MR}, 
\cite{Mon}, \cite{Moo1}, \cite{Moo2}, \cite{PS} and \cite{W}.
Details of Coxeter groups and Coxeter systems are 
found in \cite{Bo}, \cite{Br} and \cite{Hu}, and 
details of Davis complexes which are CAT(0) spaces defined by Coxeter systems 
and their boundaries are 
found in \cite{D1}, \cite{D2} and \cite{M}.
We can find some recent research on boundaries of Coxeter groups 
in \cite{CF}, \cite{Dr}, \cite{Dr2}, \cite{F}, \cite{Hos}, \cite{MRT}.
Every cocompact discrete reflection group of a geodesic space 
becomes a Coxeter group (cf.\ \cite{Ho5}), and 
we say that a Coxeter group $W$ acts geometrically on a CAT(0) space $X$ by reflections 
if the Coxeter group $W$ is a reflection group of $X$ (cf.\ \cite{Ho6}).

% Also CAT(0) cubical complex 'ɂ'¢'Ä \cite{}
% splitting theorem
% rank-one 
% JSJ decomposition 
%
% Coxeter group acts CAT(0) spaces $X$ and $Y$ by reflections?
% and Davis comlex?
% 
% cube complex case?
% cf. open problems in web
%

%%%%%%%%%%%%%%%%%%%%%%%%%%%%%%%%%%%%%
%             REFERENCES            %
%%%%%%%%%%%%%%%%%%%%%%%%%%%%%%%%%%%%%
%

%
\end{document}